\documentclass[10pt]{article}
\usepackage{amsmath}
\usepackage{amssymb}
\usepackage{amsthm}
\usepackage{mathrsfs}
\usepackage[compress,noadjust]{cite}
\numberwithin{equation}{section}

\begin{document}
\title{Singular and fractional integral operators with variable kernels on the weak Hardy spaces}
\author{Hua Wang \footnote{E-mail address: wanghua@pku.edu.cn.}\\
\footnotesize{Department of Mathematics, Zhejiang University, Hangzhou 310027, P. R. China}}
\date{}
\maketitle

\begin{abstract}
In this paper, by using the decomposition theorem for weak Hardy spaces, we will obtain the boundedness properties of some integral operators with variable kernels on these spaces, under some Dini type conditions imposed on the variable kernel $\Omega(x,z)$. 
\end{abstract}

\section{Introduction}

Let $S^{n-1}$ be the unit sphere in $\mathbb R^n$($n\ge2$) equipped with the normalized Lebesgue measure $d\sigma$. A function $\Omega(x,z)$ defined on $\mathbb R^n\times\mathbb R^n$ is said to belong to $L^{\infty}(\mathbb R^n)\times L^r(S^{n-1})$, $r\ge1$, if it satisfies the following conditions:

(1) for all $\lambda>0$ and $x,z\in\mathbb R^n$, $\Omega(x,\lambda z)=\Omega(x,z)$;

(2) $\big\|\Omega\big\|_{L^\infty(\mathbb R^n)\times L^r(S^{n-1})}:=\sup_{x\in\mathbb R^n}\left(\int_{S^{n-1}}|\Omega(x,z')|^r\,d\sigma(z')\right)^{1/r}<\infty$;

(3) for any $x\in\mathbb R^n$, $\int_{S^{n-1}}\Omega(x,z')\,d\sigma(z')=0$,

\noindent where $z'=z/{|z|}$ for any $z\in\mathbb R^n\backslash\{0\}$. Set $K(x,z)=\frac{\Omega(x,z')}{|z|^n}$. In this paper, we consider the singular integral operator with variable kernel which is defined by
\begin{equation}
T_{\Omega}f(x)=\mbox{\upshape{P.V.}}\int_{\mathbb R^n}K(x,x-y)f(y)\,dy.
\end{equation}

In 1955, Calder\'on and Zygmund \cite{cal1,cal2} investigated the $L^2$ boundedness of singular integral operators with variable kernels. They found that these operators $T_\Omega$ are closely related to the problem about second order elliptic partial differential equations with variable coefficients.  In \cite{cal1}, Calder\'on and Zygmund proved the following theorem (see also \cite{cal3}).

\newtheorem*{thma}{Theorem A}

\begin{thma}
Suppose that $\Omega(x,z)\in L^\infty(\mathbb R^n)\times L^r(S^{n-1})$ with $r>{2(n-1)}/n$, and satisfies $(1)$--$(3)$. Then there exists a constant $C>0$ independent of $f$ such that
\begin{equation*}
\big\|T_{\Omega}(f)\big\|_{L^2}\le C\big\|f\big\|_{L^2}.
\end{equation*}
\end{thma}

For $0<\alpha<n$ and $\Omega(x,z)\in L^{\infty}(\mathbb R^n)\times L^r(S^{n-1})$ with $r\ge1$, we set $K_{\alpha}(x,z)=\frac{\Omega(x,z')}{|z|^{n-\alpha}}$. Then the fractional integral operator with variable kernel is defined as follows:
\begin{equation}
T_{\Omega,\alpha}f(x)=\int_{\mathbb R^n}K_{\alpha}(x,x-y)f(y)\,dy.
\end{equation}
In 1971, Muckenhoupt and Wheeden \cite{muckenhoupt} studied the $L^p$--$L^q$ boundedness of $T_{\Omega,\alpha}$ when $0<\alpha<n$, and obtained the following result (here, and in what follows we shall denote the conjugate exponent of $p>1$ by $p'=p/{(p-1)}$):

\newtheorem*{thmb}{Theorem B}

\begin{thmb}
Let $0<\alpha<n$, $1<p<n/{\alpha}$ and $1/q=1/p-\alpha/n$. Suppose that $\Omega(x,z)\in L^\infty(\mathbb R^n)\times L^r(S^{n-1})$ with $r>p'$, and satisfies $(1)$--$(2)$. Then there exists a constant $C>0$ independent of $f$ such that
\begin{equation*}
\big\|T_{\Omega,\alpha}(f)\big\|_{L^q}\le C\big\|f\big\|_{L^p}.
\end{equation*}
\end{thmb}

On the other hand, the weak $H^p$ spaces have first appeared in the work of Fefferman, Rivi\`ere and Sagher \cite{cfefferman}, which are the intermediate spaces between two Hardy spaces through the real method of interpolation. The atomic decomposition characterization of weak $H^1$ space on $\mathbb R^n$ was given by Fefferman and Soria in \cite{rfefferman}. Later, Liu \cite{liu1} established the weak $H^p$ spaces on homogeneous groups for the whole range $0<p\le1$. The corresponding results related to $\mathbb R^n$ can be found in \cite{lu}. For the continuity properties of some kinds of operators on weak Hardy spaces, we refer the readers to \cite{ding1,ding2,ding3,ding4,ding5,liu2,tao}.

In \cite{ding1}, the authors considered the boundedness of $T_\Omega$ and $T_{\Omega,\alpha}$ on the weak Hardy spaces $WH^1(\mathbb R^n)$, under certain smoothness conditions on the variable kernel $\Omega(x,z)$. Motivated by \cite{ding1}, the main purpose of this paper is to establish the boundedness properties of $T_\Omega$ and $T_{\Omega,\alpha}$ on the spaces $WH^p(\mathbb R^n)$, under the assumptions that $\Omega(x,z)$ satisfy some Dini type conditions (see Section 2 for its definition). We now formulate our main results as follows.

\newtheorem{theorem}{Theorem}[section]

\begin{theorem}
Let $0<\beta\le 1$ and $n/{(n+\beta)}<p\le1$. Suppose that $\Omega\in\mbox{Din}^r_\beta(S^{n-1})$ with $r>{2(n-1)}/n$, then there exists a constant $C>0$ independent of $f$ such that
\begin{equation*}
\big\|T_\Omega(f)\big\|_{WL^p}\le C\big\|f\big\|_{WH^p}.
\end{equation*}
\end{theorem}

\begin{theorem}
Let $0<\alpha\le 1$, $n/{(n+\alpha)}<p\le1$ and $1/q=1/p-\alpha/n$. Suppose that $\Omega\in\mbox{Din}^q_\alpha(S^{n-1})$, then there exists a constant $C>0$ independent of $f$ such that
\begin{equation*}
\big\|T_{\Omega,\alpha}(f)\big\|_{WL^q}\le C\big\|f\big\|_{WH^p}.
\end{equation*}
\end{theorem}

\begin{theorem}
Let $0<\alpha<\beta\le 1$, $n/{(n+\beta)}<p\le n/{(n+\alpha)}$ and $1/q=1/p-\alpha/n$. Suppose that $\Omega\in\mbox{Din}^r_\beta(S^{n-1})$ with $r>n/{(n-\alpha)}$, then there exists a constant $C>0$ independent of $f$ such that
\begin{equation*}
\big\|T_{\Omega,\alpha}(f)\big\|_{WL^q}\le C\big\|f\big\|_{WH^p}.
\end{equation*}
\end{theorem}

\section{Notations and preliminaries}

For $0<p<\infty$, we denote by $L^p(\mathbb R^n)$ the classical Lebesgue spaces of all functions $f$ satisfying
\begin{equation}
\big\|f\big\|_{L^p}=\left(\int_{\mathbb R^n}|f(x)|^p\,dx\right)^{1/p}<\infty.
\end{equation}
When $p=\infty$, $L^\infty(\mathbb R^n)$ will be defined as follows:
\begin{equation}
\big\|f\big\|_{L^\infty}=\underset{x\in\mathbb R^n}{\mbox{ess\,sup}}\,|f(x)|<\infty.
\end{equation}
We also denote by $WL^p(\mathbb R^n)$ the weak $L^p$ spaces consisting of all measurable functions $f$ such that
\begin{equation}
\big\|f\big\|_{WL^p}=\sup_{\lambda>0}\lambda\cdot \big|\big\{x\in\mathbb R^n:|f(x)|>\lambda \big\}\big|^{1/p}<\infty.
\end{equation}

Let us first recall the definitions of the integral Dini type conditions. In \cite{ding6,ding7}, Ding et al. introduced some definitions about the variable kernel $\Omega(x,z)$ when they studied the $H^1$--$L^1$ boundedness of Marcinkiewicz integrals. Replacing the condition (2) mentioned before, they strengthened it to the condition

$(2')$\, $\sup_{x\in\mathbb R^n\atop\rho\ge0}\left(\int_{S^{n-1}}|\Omega(x+\rho z',z')|^r\,d\sigma(z')\right)^{1/r}<\infty$.

For $\Omega(x,z)\in L^{\infty}(\mathbb R^n)\times L^r(S^{n-1})$ and $r\ge1$, a function $\Omega(x,z)$ is said to satisfy the $L^r$-Dini condition if the conditions $(1)$, $(2')$, $(3)$ hold and
\begin{equation}
\int_0^1\frac{\omega_r(\delta)}{\delta}\,d\delta<\infty,
\end{equation}
where $\omega_r(\delta)$ is the integral modulus of continuity of order $r$ of $\Omega$ defined by
\begin{equation}
\omega_r(\delta):=\sup_{x\in\mathbb R^n\atop\rho\ge0}\bigg(\int_{S^{n-1}}\sup_{y'\in S^{n-1}\atop |y'-z'|\le\delta}\big|\Omega(x+\rho z',y')-\Omega(x+\rho z',z')\big|^rd\sigma(z')\bigg)^{1/r}.
\end{equation}
In order to obtain the $H^p$--$L^p$ boundedness of $T_{\Omega}$, Lee et al. \cite{lee} generalized the $L^r$-Dini condition by replacing (2.4) to the following stronger condition (see also \cite{lin})
\begin{equation}
\int_0^1\frac{\omega_r(\delta)}{\delta^{1+\alpha}}\,d\delta<\infty,\quad 0\le\alpha\le1.
\end{equation}
If $\Omega$ satisfies (2.6) for some $r\ge1$ and $0\le\alpha\le1$, we say that it satisfies the $L^{r,\alpha}$-Dini condition. For the special case $\alpha=0$, it reduces to the $L^r$-Dini condition. For $0\le\beta<\alpha\le1$, if $\Omega$ satisfies the $L^{r,\alpha}$-Dini condition, then it also satisfies the $L^{r,\beta}$-Dini condition. We thus denote by $\mbox{Din}^r_\alpha(S^{n-1})$ the class of all functions which satisfy the $L^{r,\beta}$-Dini condition for all $0<\beta<\alpha$.
Following the same arguments as in the proof of Lemma 5 in \cite{kurtz}, we can also establish the following lemma on the variable kernel $\Omega(x,z)$ (See \cite{ding7} and \cite{lee}).

\newtheorem{lemma}[theorem]{Lemma}

\begin{lemma}
Let $0\le\alpha<n$ and $r\ge1$. Suppose that $\Omega(x,z)\in L^\infty(\mathbb R^n)\times L^r(S^{n-1})$ satisfies the $L^r$-Dini condition of this section. If there exists a constant $0<\gamma\le 1/2$ such that $|y|<\gamma R$, then for any $x_0\in\mathbb R^n$, we have
\begin{equation*}
\begin{split}
\bigg(&\int_{R\le|x|<2R}\big|K_{\alpha}(x+x_0,x-y)-K_{\alpha}(x+x_0,x)\big|^rdx\bigg)^{1/r}\\
&\le C\cdot R^{n/r-(n-\alpha)}\bigg(\frac{|y|}{R}+\int_{|y|/{2R}}^{|y|/R}\frac{\omega_r(\delta)}{\delta}d\delta\bigg),
\end{split}
\end{equation*}
where the constant $C>0$ is independent of $R$ and $y$. We simply denote $K_{\alpha}(x,z)$ by $K(x,z)$ when $\alpha=0$.
\end{lemma}

Now let us turn to the weak Hardy spaces. We write $\mathscr S(\mathbb R^n)$ to denote the Schwartz space of all rapidly decreasing infinitely differentiable functions and $\mathscr S'(\mathbb R^n)$ to denote the space of all tempered distributions, i.e., the topological dual of $\mathscr S(\mathbb R^n)$. Let $0<p\le1$ and $N=[n(1/p-1)]$. Define
\begin{equation*}
\mathscr A_{N}=\Big\{\varphi\in\mathscr S(\mathbb R^n):\sup_{x\in\mathbb R^n}\sup_{|\alpha|\le N+1}(1+|x|)^{N+n+1}\big|D^\alpha\varphi(x)\big|\le1\Big\},
\end{equation*}
where $\alpha=(\alpha_1,\dots,\alpha_n)\in(\mathbb N\cup\{0\})^n$, $|\alpha|=\alpha_1+\dots+\alpha_n$, and
\begin{equation*}
D^\alpha\varphi=\frac{\partial^{|\alpha|}\varphi}{\partial x^{\alpha_1}_1\cdots\partial x^{\alpha_n}_n}.
\end{equation*}
For any given $f\in\mathscr S'(\mathbb R^n)$, the grand maximal function of $f$ is defined by
\begin{equation*}
G f(x)=\sup_{\varphi\in\mathscr A_{N}}\sup_{|y-x|<t}\big|(\varphi_t*f)(y)\big|.
\end{equation*}
Then we can define the weak Hardy space $WH^p(\mathbb R^n)$ by $WH^p(\mathbb R^n)=\big\{f\in\mathscr S'(\mathbb R^n):G(f)\in WL^p(\mathbb R^n)\big\}$. Moreover, we set $\big\|f\big\|_{WH^p}=\big\|G(f)\big\|_{WL^p}$.

We need the following atomic decomposition theorem for weak Hardy spaces $WH^p(\mathbb R^n)$ given in \cite{liu1}
(see also \cite{lu}).

\begin{theorem}
Let $0<p\le1$. For every $f\in WH^p(\mathbb R^n)$, then there exists a sequence of bounded measurable functions $\{f_k\}_{k=-\infty}^\infty$ such that

$(i)$ $f=\sum_{k=-\infty}^\infty f_k$ in the sense of distributions.

$(ii)$ Each $f_k$ can be further decomposed into $f_k=\sum_i b^k_i$, where $\{b^k_i\}$ satisfies

\quad $(a)$ Each $b^k_i$ is supported in a cube $Q^k_i$ with $\sum_{i}\big|Q^k_i\big|\le c2^{-kp}$, and $\sum_i\chi_{Q^k_i}(x)\le c$. Here $\chi_E$ denotes the characteristic function of the set $E$ and $c\sim\big\|f\big\|_{WH^p}^p;$

\quad $(b)$ $\big\|b^k_i\big\|_{L^\infty}\le C2^k$, where $C>0$ is independent of $i$ and $k\,;$

\quad $(c)$ $\int_{\mathbb R^n}b^k_i(x)x^\gamma\,dx=0$ for every multi-index $\gamma$ with $|\gamma|\le[n(1/p-1)]$.

Conversely, if $f\in\mathscr S'(\mathbb R^n)$ has a decomposition satisfying $(i)$ and $(ii)$, then $f\in WH^p(\mathbb R^n)$. Moreover, we have $\big\|f\big\|_{WH^p}^p\sim c.$
\end{theorem}

Throughout this article $C$ always denotes a positive constant, which is independent of the main parameters and not necessarily the same at each occurrence.

\section{Proof of Theorem 1.1}

\begin{proof}[Proof of Theorem 1.1]
For any given $\lambda>0$, we may choose $k_0\in\mathbb Z$ such that $2^{k_0}\le\lambda<2^{k_0+1}$. For every $f\in WH^p(\mathbb R^n)$, then by Theorem 2.2, we can write
\begin{equation*}
f=\sum_{k=-\infty}^\infty f_k=\sum_{k=-\infty}^{k_0} f_k+\sum_{k=k_0+1}^\infty f_k:=F_1+F_2,
\end{equation*}
where $F_1=\sum_{k=-\infty}^{k_0} f_k=\sum_{k=-\infty}^{k_0}\sum_i b^k_i$, $F_2=\sum_{k=k_0+1}^\infty f_k=\sum_{k=k_0+1}^\infty\sum_i b^k_i$ and $\{b^k_i\}$ satisfies $(a)$--$(c)$ in Theorem 2.2. Then we have \begin{equation*}
\begin{split}
&\lambda^p\cdot \big|\big\{x\in\mathbb R^n:\big|T_{\Omega}(f)(x)\big|>\lambda\big\}\big|\\
\le\,&\lambda^p\cdot \big|\big\{x\in\mathbb R^n:\big|T_{\Omega}(F_1)(x)\big|>\lambda/2\big\}\big|
+\lambda^p\cdot \big|\big\{x\in\mathbb R^n:\big|T_{\Omega}(F_2)(x)\big|>\lambda/2\big\}\big|\\
=\,&I_1+I_2.
\end{split}
\end{equation*}
First we claim that the following inequality holds:
\begin{equation}
\big\|F_1\big\|_{L^2}\le C\cdot\lambda^{1-p/2}\big\|f\big\|^{p/2}_{WH^p}.
\end{equation}
In fact, since supp\,$b^k_i\subseteq Q^k_i=Q\big(x^k_i,r^k_i\big)$ and $\big\|b^k_i\big\|_{L^\infty}\le C 2^k$ by Theorem 2.2, where $Q\big(x^k_i,r^k_i\big)$ denotes the cube centered at $x^k_i$ with side length $r^k_i$ and all cubes are assumed to have their sides parallel to the coordinate axes. Hence, it follows from Minkowski's integral inequality that
\begin{equation*}
\begin{split}
\big\|F_1\big\|_{L^2}&\le\sum_{k=-\infty}^{k_0}\sum_i\big\|b^k_i\big\|_{L^2}\\
&\le\sum_{k=-\infty}^{k_0}\sum_i\big\|b^k_i\big\|_{L^\infty}\big|Q^k_i\big|^{1/2}.
\end{split}
\end{equation*}
For each $k\in\mathbb Z$, by using the bounded overlapping property of the cubes $\{Q^k_i\}$ and the fact that $1-p/2>0$, we thus obtain
\begin{equation*}
\begin{split}
\big\|F_1\big\|_{L^2}&\le C\sum_{k=-\infty}^{k_0}2^k\Big(\sum_i \big|Q^k_i\big|\Big)^{1/2}\\
&\le C\sum_{k=-\infty}^{k_0}2^{k(1-p/2)}\big\|f\big\|^{p/2}_{WH^p}\\
&\le C\sum_{k=-\infty}^{k_0}2^{(k-k_0)(1-p/2)}\cdot\lambda^{1-p/2}\big\|f\big\|^{p/2}_{WH^p}\\
&\le C\cdot\lambda^{1-p/2}\big\|f\big\|^{p/2}_{WH^p}.
\end{split}
\end{equation*}
Note that $\Omega\in\mbox{Din}^r_\beta(S^{n-1})$ with $r>{2(n-1)}/n$, then we know that $T_{\Omega}$ is bounded on $L^2(\mathbb R^n)$ according to Theorem A. This fact together with Chebyshev's inequality and the inequality (3.1) yields
\begin{align}
I_1&\le \lambda^p\cdot\frac{4}{\lambda^2}\big\|T_{\Omega}(F_1)\big\|^2_{L^2}\notag\\
&\le C\cdot\lambda^{p-2}\big\|F_1\big\|^2_{L^2}\notag\\
&\le C\big\|f\big\|^{p}_{WH^p}.
\end{align}
We now turn our attention to the estimate of $I_2$. Setting
\begin{equation*}
A_{k_0}=\bigcup_{k=k_0+1}^\infty\bigcup_i \widetilde{Q^k_i},
\end{equation*}
where $\widetilde{Q^k_i}=Q\big(x^k_i,\tau^{{p(k-k_0)}/n}(2\sqrt n)r^k_i\big)$ and $\tau$ is a fixed positive number such that $1<\tau<2$. Thus, we can further decompose $I_2$ as
\begin{equation*}
\begin{split}
I_2&\le\lambda^p\cdot \big|\big\{x\in A_{k_0}:|T_{\Omega}(F_2)(x)|>\lambda/2\big\}\big|+
\lambda^p\cdot \big|\big\{x\in (A_{k_0})^c:|T_{\Omega}(F_2)(x)|>\lambda/2\big\}\big|\\
&=I'_2+I''_2.
\end{split}
\end{equation*}
For the term $I'_2$, we can deduce that
\begin{align}
I'_2&\le\lambda^p\sum_{k=k_0+1}^\infty\sum_i \big|\widetilde{Q^k_i}\big|\notag\\
&\le C\cdot\lambda^p\sum_{k=k_0+1}^\infty\tau^{p(k-k_0)}\sum_i \big|Q^k_i\big|\notag\\
&\le C\big\|f\big\|^{p}_{WH^p}\sum_{k=k_0+1}^\infty\Big(\frac{\tau}{2}\Big)^{p(k-k_0)}\notag\\
&\le C\big\|f\big\|^{p}_{WH^p}.
\end{align}
On the other hand, it follows directly from Chebyshev's inequality that
\begin{equation*}
\begin{split}
I''_2&\le 2^p\int_{(A_{k_0})^c}\big|T_{\Omega}(F_2)(x)\big|^p\,dx\\
&\le 2^p
\sum_{k=k_0+1}^\infty\sum_i\int_{\big(\widetilde{Q^k_i}\big)^c}\big|T_{\Omega}\big(b^k_i\big)(x)\big|^p\,dx\\
&= 2^p
\sum_{k=k_0+1}^\infty\sum_i {\mathcal J}^k_i.
\end{split}
\end{equation*}
Now denote $\tau^k_{i,\ell}=2^{\ell-1}\tau^{{p(k-k_0)}/n}\sqrt{n}r^k_i$ and
$$E^k_{i,\ell}=\big\{x\in\mathbb R^n:\tau^k_{i,\ell}\le|x-x^k_i|<2\tau^k_{i,\ell}\big\},\quad \ell=1,2,\ldots.$$
An application of H\"older's inequality gives us that
\begin{equation*}
\begin{split}
{\mathcal J}^k_i&\le\sum_{\ell=1}^\infty\int_{E^k_{i,\ell}}\big|T_{\Omega}\big(b^k_i\big)(x)\big|^p\,dx\\
&\le\sum_{\ell=1}^\infty\Big|E^k_{i,\ell}\Big|^{1-p}
\bigg(\int_{E^k_{i,\ell}}\big|T_{\Omega}\big(b^k_i\big)(x)\big|\,dx\bigg)^p.
\end{split}
\end{equation*}
Observe that $[n(1/p-1)]=0$ by our assumptions. Thus, by the cancellation condition of $b^k_i\in L^\infty(\mathbb R^n)$, we get
\begin{equation*}
\begin{split}
\int_{E^k_{i,\ell}}\big|T_{\Omega}\big(b^k_i\big)(x)\big|\,dx
&=\int_{E^k_{i,\ell}}\left|\int_{Q^k_i}\Big[K\big(x,x-y\big)-K\big(x,x-x^k_i\big)\Big]b^k_i(y)\,dy\right|\,dx\\
&\le\int_{Q^k_i}\bigg\{\int_{E^k_{i,\ell}}\Big|K\big(x,x-y\big)-K\big(x,x-x^k_i\big)\Big|\,dx\bigg\}
\big|b^k_i(y)\big|\,dy\\
&\le\big\|b^k_i\big\|_{L^\infty}\int_{Q^k_i}\bigg\{\int_{E^k_{i,\ell}}
\Big|K\big(x,x-y\big)-K\big(x,x-x^k_i\big)\Big|\,dx\bigg\}\,dy.
\end{split}
\end{equation*}
When $y\in Q^k_i$ and $x\in\big(\widetilde{Q^k_i}\big)^c$, then a trivial computation shows that for all $i$ and $k$,
\begin{equation}
\big|x-x^k_i\big|\ge\tau^{{(k-k_0)}/{(n+\alpha)}}\sqrt{n} r^k_i>\sqrt{n} r^k_i\ge 2\big|y-x^k_i\big|.
\end{equation}
Using H\"older's inequality, the estimate (3.4) and Lemma 2.1, we can see that for any $y\in Q^k_i$, the integral in the brace of the above expression is dominated by
\begin{align}
&\bigg(\int_{E^k_{i,\ell}}\Big|K\big(x,x-y\big)-K\big(x,x-x^k_i\big)\Big|^r\,dx\bigg)^{1/r}
\bigg(\int_{E^k_{i,\ell}}1\,dx\bigg)^{1/{r'}}\notag\\
\le\, &C\cdot\Big|E^k_{i,\ell}\Big|^{1/{r'}}
\left(\int_{\tau^k_{i,\ell}\le|x|<2\tau^k_{i,\ell}}\Big|K\big(x+x^k_i,x-(y-x^k_i)\big)
-K\big(x+x^k_i,x\big)\Big|^r\,dx\right)^{1/r}\notag\\
\le\, &C\cdot\Big|E^k_{i,\ell}\Big|^{1/{r'}}
\cdot\Big(\tau^k_{i,\ell}\Big)^{-n/{r'}}
\left(\frac{|y-x^k_i|}{\tau^k_{i,\ell}}+
\int_{|y-x^k_i|/{2\tau^k_{i,\ell}}}^{|y-x^k_i|/{\tau^k_{i,\ell}}}\frac{\omega_r(\delta)}{\delta}\,d\delta\right)
\notag\\
\le\, &C\cdot\Big|E^k_{i,\ell}\Big|^{1/{r'}}
\cdot\Big(\tau^k_{i,\ell}\Big)^{-n/{r'}}
\left(\frac{|y-x^k_i|}{\tau^k_{i,\ell}}+\frac{|y-x^k_i|^\beta}{(\tau^k_{i,\ell})^\beta}\times
\int_{|y-x^k_i|/{2\tau^k_{i,\ell}}}^{|y-x^k_i|/{\tau^k_{i,\ell}}}
\frac{\omega_r(\delta)}{\delta^{1+\beta}}\,d\delta\right)
\notag\\
\le\, & C\cdot\Big(2\tau^k_{i,\ell}\Big)^{n/{r'}}
\cdot\Big(\tau^k_{i,\ell}\Big)^{-n/{r'}}
\left(\frac{1}{2^{\ell}\tau^{{p(k-k_0)}/n}}+
\Big[\frac{1}{2^{\ell}\tau^{{p(k-k_0)}/n}}\Big]^\beta
\int_{0}^{1}\frac{\omega_r(\delta)}{\delta^{1+\beta}}\,d\delta\right)
\notag\\
\le\, & C\cdot\left(1+\int_{0}^{1}\frac{\omega_r(\delta)}{\delta^{1+\beta}}\,d\delta\right)
\times\left(\frac{1}{2^{\ell}\tau^{{p(k-k_0)}/n}}\right)^\beta.
\end{align}
Recall that $\big\|b^k_i\big\|_{L^\infty}\le C 2^k$. From the above estimate (3.5), it follows that for all $i$ and $k$,
\begin{equation*}
\begin{split}
{\mathcal J}^k_i&\le C\cdot2^{kp}\sum_{\ell=1}^\infty \Big|E^k_{i,\ell}\Big|^{1-p}\cdot\Big|Q^k_i\Big|^p
\left(\frac{1}{2^{\ell}\tau^{{p(k-k_0)}/n}}\right)^{\beta p}\\
&\le C\cdot2^{kp}\sum_{\ell=1}^\infty \Big|Q^k_i\Big|^{1-p}\cdot\Big|Q^k_i\Big|^p
\Big(2^{\ell}\tau^{{p(k-k_0)}/n}\Big)^{n(1-p)-\beta p}\\
&\le C\cdot2^{kp}\cdot\big|Q^k_i\big|\Big(\tau^{{p(k-k_0)}/n}\Big)^{n(1-p)-\beta p},
\end{split}
\end{equation*}
where the last inequality holds since $p>n/{(n+\beta)}$. Therefore
\begin{align}
I''_2&\le C\sum_{k=k_0+1}^\infty\sum_i 2^{kp}\cdot\big|Q^k_i\big|
\Big(\tau^{{p(k-k_0)}/n}\Big)^{n-(n+\beta)p}\notag\\
&\le C\big\|f\big\|^{p}_{WH^p}\sum_{k=k_0+1}^\infty\Big(\tau^{{p(k-k_0)}/n}\Big)^{n-(n+\beta)p}\notag\\
&\le C\big\|f\big\|^{p}_{WH^p}\sum_{k=1}^\infty\Big(\tau^{{pk}/n}\Big)^{n-(n+\beta)p}\notag\\
&\le C\big\|f\big\|^{p}_{WH^p}.
\end{align}
Combining the above inequality (3.6) with (3.2)--(3.3) and taking the supremum over all $\lambda>0$, and then taking $p$-th root on both sides, we complete the proof of Theorem 1.1.
\end{proof}

\section{Proof of Theorem 1.2}

\begin{proof}[Proof of Theorem 1.2]
For any fixed $\lambda>0$, we may choose $k_0\in\mathbb Z$ satisfying $2^{k_0}\le\xi<2^{k_0+1}$, where we define $\xi=\lambda^{q/p}\big\|f\big\|^{1-q/p}_{WH^p}$. For every $f\in WH^p(\mathbb R^n)$, then in view of Theorem 2.2, we can write
\begin{equation*}
f=\sum_{k=-\infty}^\infty f_k=\sum_{k=-\infty}^{k_0} f_k+\sum_{k=k_0+1}^\infty f_k:=F_1+F_2,
\end{equation*}
where $F_1=\sum_{k=-\infty}^{k_0} f_k=\sum_{k=-\infty}^{k_0}\sum_i b^k_i$, $F_2=\sum_{k=k_0+1}^\infty f_k=\sum_{k=k_0+1}^\infty\sum_i b^k_i$ and $\{b^k_i\}$ satisfies $(a)$--$(c)$ in Theorem 2.2. Then we have
\begin{equation*}
\begin{split}
&\lambda^q\cdot \big|\big\{x\in\mathbb R^n:\big|T_{\Omega,\alpha}(f)(x)\big|>\lambda\big\}\big|\\
\le\,&\lambda^q\cdot \big|\big\{x\in\mathbb R^n:\big|T_{\Omega,\alpha}(F_1)(x)\big|>\lambda/2\big\}\big|
+\lambda^q\cdot \big|\big\{x\in\mathbb R^n:\big|T_{\Omega,\alpha}(F_2)(x)\big|>\lambda/2\big\}\big|\\
=\,&J_1+J_2.
\end{split}
\end{equation*}
If $0<\alpha\le 1$, $n/{(n+\alpha)}<p\le1$ and $1/q=1/p-\alpha/n$, then $q>1$. Thus, we are able to choose $p_1$ such that $1<p_1<n/{\alpha}$ and $q>p_1'>1$. Then we take $q_1>p_1>1$ such that $1/{q_1}=1/{p_1}-\alpha/n$. Similar to the proof of Theorem 1.1, we first claim that the following inequality holds:
\begin{equation}
\big\|F_1\big\|_{L^{p_1}}\le C\cdot\xi^{1-p/{p_1}}\big\|f\big\|^{p/{p_1}}_{WH^p}.
\end{equation}
Indeed, since supp\,$b^k_i\subseteq Q^k_i=Q\big(x^k_i,r^k_i\big)$ and $\big\|b^k_i\big\|_{L^\infty}\le C 2^k$ by Theorem 2.2, then by using Minkowski's integral inequality, we get
\begin{equation*}
\begin{split}
\big\|F_1\big\|_{L^{p_1}}&\le\sum_{k=-\infty}^{k_0}\sum_i\big\|b^k_i\big\|_{L^{p_1}}\\
&\le\sum_{k=-\infty}^{k_0}\sum_i\big\|b^k_i\big\|_{L^\infty}\big|Q^k_i\big|^{1/{p_1}}.
\end{split}
\end{equation*}
For each $k\in\mathbb Z$, by using the finitely overlapping property of the cubes $\{Q^k_i\}$ and the fact that $1-p/{p_1}>0$, we thus obtain
\begin{equation*}
\begin{split}
\big\|F_1\big\|_{L^{p_1}}&\le C\sum_{k=-\infty}^{k_0}2^k\Big(\sum_i \big|Q^k_i\big|\Big)^{1/{p_1}}\\
&\le C\sum_{k=-\infty}^{k_0}2^{k(1-p/{p_1})}\big\|f\big\|^{p/{p_1}}_{WH^p}\\
&\le C\sum_{k=-\infty}^{k_0}2^{(k-k_0)(1-p/{p_1})}\cdot\xi^{1-p/{p_1}}\big\|f\big\|^{p/{p_1}}_{WH^p}\\
&\le C\cdot\xi^{1-p/{p_1}}\big\|f\big\|^{p/{p_1}}_{WH^p}.
\end{split}
\end{equation*}
Notice that $\Omega\in\mbox{Din}^q_\alpha(S^{n-1})$ with $q>p_1'$, then we know that $T_{\Omega,\alpha}$ is bounded from $L^{p_1}(\mathbb R^n)$ to $L^{q_1}(\mathbb R^n)$ according to Theorem B. This fact along with Chebyshev's inequality and the inequality (4.1) implies
\begin{align}
J_1&\le\lambda^q\cdot\Big(\frac{2}{\lambda}\Big)^{q_1}
\big\|T_{\Omega,\alpha}(F_1)\big\|^{q_1}_{L^{q_1}}\notag\\
&\le C\cdot\lambda^{q-q_1}\big\|F_1\big\|^{q_1}_{L^{p_1}}\notag\\
&\le C\cdot\lambda^{q-q_1}\xi^{(1-p/{p_1})q_1}\big\|f\big\|^{{pq_1}/{p_1}}_{WH^p}\notag\\
&\le C\cdot\lambda^{q-q_1}\Big(\lambda^{q/p}\big\|f\big\|^{1-q/p}_{WH^p}\Big)^{(1-p/{p_1})q_1}
\big\|f\big\|^{{pq_1}/{p_1}}_{WH^p}.
\end{align}
Note that $1/p-1/q=1/{p_1}-1/{q_1}=\alpha/n$, then it is easy to check that
\begin{equation*}
\begin{split}
q-q_1+q/p\cdot(1-p/{p_1})q_1&=q-q_1+qq_1\cdot(1/p-1/{p_1})\\
&=q-q_1+qq_1\cdot(1/q-1/{q_1})\\
&=0
\end{split}
\end{equation*}
and
\begin{equation*}
\begin{split}
(1-q/p)\cdot(1-p/{p_1})q_1+{pq_1}/{p_1}&=(1-q/p)\cdot(1-p/{p_1})q_1-(1-p/{p_1})q_1+q_1\\
&=q/p\cdot(p/{p_1}-1)q_1+q_1\\
&=qq_1\cdot(1/{q_1}-1/q)+q_1\\
&=q.
\end{split}
\end{equation*}
Hence, by the inequality (4.2), we have
\begin{equation}
J_1\le C\big\|f\big\|^{q}_{WH^p}.
\end{equation}
Let us now turn our attention to the estimate of $J_2$. Setting
\begin{equation*}
A_{k_0}=\bigcup_{k=k_0+1}^\infty\bigcup_i \widetilde{Q^k_i},
\end{equation*}
where $\widetilde{Q^k_i}=Q\big(x^k_i,\tau^{{p(k-k_0)}/n}(2\sqrt n)r^k_i\big)$ and $\tau$ is also a fixed positive number such that $1<\tau<2$. Thus, we can further split $J_2$ into
\begin{equation*}
\begin{split}
J_2&\le\lambda^q\cdot \big|\big\{x\in A_{k_0}:|T_{\Omega,\alpha}(F_2)(x)|>\lambda/2\big\}\big|+
\lambda^q\cdot \big|\big\{x\in (A_{k_0})^c:|T_{\Omega,\alpha}(F_2)(x)|>\lambda/2\big\}\big|\\
&=J'_2+J''_2.
\end{split}
\end{equation*}
For the term $J'_2$, we can see that
\begin{align}
J'_2&\le\lambda^q\sum_{k=k_0+1}^\infty\sum_i \big|\widetilde{Q^k_i}\big|\notag\\
&\le C\cdot\lambda^q\sum_{k=k_0+1}^\infty\tau^{p(k-k_0)}\sum_i \big|Q^k_i\big|\notag\\
&\le C\cdot\lambda^q\cdot\xi^{-p}
\big\|f\big\|^{p}_{WH^p}\sum_{k=k_0+1}^\infty\Big(\frac{\tau}{2}\Big)^{p(k-k_0)}\notag\\
&\le C\big\|f\big\|^{q}_{WH^p}.
\end{align}
For the term $J''_2$, we denote $\tau^k_{i,\ell}=2^{\ell-1}\tau^{{p(k-k_0)}/n}\sqrt{n}r^k_i$ and
$$E^k_{i,\ell}=\big\{x\in\mathbb R^n:\tau^k_{i,\ell}\le|x-x^k_i|<2\tau^k_{i,\ell}\big\},\quad \ell=1,2,\ldots.$$
Then it follows directly from Chebyshev's inequality that
\begin{equation*}
\begin{split}
J''_2&\le 2^q\int_{(A_{k_0})^c}\big|T_{\Omega,\alpha}(F_2)(x)\big|^q\,dx\\
&\le 2^q\sum_{k=k_0+1}^\infty\sum_i
\int_{\big(\widetilde{Q^k_i}\big)^c}\big|T_{\Omega,\alpha}\big(b^k_i\big)(x)\big|^q\,dx\\
&\le 2^q\sum_{k=k_0+1}^\infty\sum_i
\sum_{\ell=1}^\infty\int_{E^k_{i,\ell}}\big|T_{\Omega,\alpha}\big(b^k_i\big)(x)\big|^q\,dx.
\end{split}
\end{equation*}
Observe that $[n(1/p-1)]=0$ and $q>1$. Hence, by using H\"older's inequality with exponent $q$ and the cancellation condition of $b^k_i\in L^\infty(\mathbb R^n)$, we deduce that
\begin{equation*}
\begin{split}
J''_2&\le C\sum_{k=k_0+1}^\infty\sum_i\sum_{\ell=1}^\infty
\int_{E^k_{i,\ell}}
\left|\int_{Q^k_i}\Big[K_{\alpha}\big(x,x-y\big)-K_{\alpha}\big(x,x-x^k_i\big)\Big]b^k_i(y)\,dy\right|^q\,dx\\
&\le C\sum_{k=k_0+1}^\infty\sum_i\sum_{\ell=1}^\infty
\int_{E^k_{i,\ell}}\left(\int_{Q^k_i}\Big|K_{\alpha}\big(x,x-y\big)-K_{\alpha}\big(x,x-x^k_i\big)\Big|^qdy\right)\\
&\times\left(\int_{Q^k_i}\big|b^k_i(y)\big|^{q'}\,dy\right)^{q/{q'}}dx\\
&\le C\sum_{k=k_0+1}^\infty\sum_i\big\|b^k_i\big\|^q_{L^\infty}\cdot\big|Q^k_i\big|^{q/{q'}}\sum_{\ell=1}^\infty
\int_{Q^k_i}\left(\int_{E^k_{i,\ell}}\Big|K_{\alpha}\big(x,x-y\big)-K_{\alpha}\big(x,x-x^k_i\big)\Big|^qdx\right)dy.
\end{split}
\end{equation*}
If $y\in Q^k_i$ and $x\in\big(\widetilde{Q^k_i}\big)^c$, then we still have $\big|x-x^k_i\big|\ge 2\big|y-x^k_i\big|$ for all $i$ and $k$ by (3.4). Since $\Omega\in\mbox{Din}^q_\alpha(S^{n-1})$ with $q>1$, then by Lemma 2.1,for any $y\in Q^k_i$, we obtain
\begin{equation*}
\begin{split}
&\left(\int_{E^k_{i,\ell}}\Big|K_{\alpha}\big(x,x-y\big)-K_{\alpha}\big(x,x-x^k_i\big)\Big|^qdx\right)^{1/q}\\
\le &\left(\int_{\tau^k_{i,\ell}\le|x|<2\tau^k_{i,\ell}}\Big|K_{\alpha}\big(x+x^k_i,x-(y-x^k_i)\big)
-K_{\alpha}\big(x+x^k_i,x\big)\Big|^q\,dx\right)^{1/q}\\
\le & C\cdot\Big(\tau^k_{i,\ell}\Big)^{n/q-(n-\alpha)}
\left(\frac{|y-x^k_i|}{\tau^k_{i,\ell}}+
\int_{|y-x^k_i|/{2\tau^k_{i,\ell}}}^{|y-x^k_i|/{\tau^k_{i,\ell}}}\frac{\omega_q(\delta)}{\delta}\,d\delta\right)\\
\le & C\cdot\Big(\tau^k_{i,\ell}\Big)^{n/q-(n-\alpha)}
\left(\frac{|y-x^k_i|}{\tau^k_{i,\ell}}+\frac{|y-x^k_i|^\alpha}{(\tau^k_{i,\ell})^\alpha}\times
\int_{|y-x^k_i|/{2\tau^k_{i,\ell}}}^{|y-x^k_i|/{\tau^k_{i,\ell}}}
\frac{\omega_q(\delta)}{\delta^{1+\alpha}}\,d\delta\right)\\
\le & C\cdot\Big(\tau^k_{i,\ell}\Big)^{n/q-(n-\alpha)}
\left(\frac{1}{2^{\ell}\tau^{{p(k-k_0)}/n}}+
\Big[\frac{1}{2^{\ell}\tau^{{p(k-k_0)}/n}}\Big]^\alpha
\int_{0}^{1}\frac{\omega_q(\delta)}{\delta^{1+\alpha}}\,d\delta\right)\\
\le & C\cdot\Big(\tau^k_{i,\ell}\Big)^{n/q-(n-\alpha)}
\left(1+\int_{0}^{1}\frac{\omega_q(\delta)}{\delta^{1+\alpha}}\,d\delta\right)
\times\left(\frac{1}{2^{\ell}\tau^{{p(k-k_0)}/n}}\right)^\alpha.
\end{split}
\end{equation*}
So we have
\begin{equation*}
\begin{split}
J''_2&\le C\sum_{k=k_0+1}^\infty\sum_i\big\|b^k_i\big\|^q_{L^\infty}\cdot\Big|Q^k_i\Big|^{q/{q'}+1}
\sum_{\ell=1}^\infty
\Big(\tau^k_{i,\ell}\Big)^{n-(n-\alpha)q}\left(\frac{1}{2^{\ell}\tau^{{p(k-k_0)}/n}}\right)^{\alpha q}\\
&\le C\sum_{k=k_0+1}^\infty2^{kq}\cdot\sum_i\Big|Q^k_i\Big|^{q}\sum_{\ell=1}^\infty
\Big|Q^k_i\Big|^{1-{(n-\alpha)q}/n}
\Big(2^{\ell}\tau^{{p(k-k_0)}/n}\Big)^{n(1-q)}\\
&\le C\sum_{k=k_0+1}^\infty2^{kq}\Big(\tau^{{p(k-k_0)}/n}\Big)^{n(1-q)}\cdot\sum_i\Big|Q^k_i\Big|^{q/p},
\end{split}
\end{equation*}
where in the last inequality we have used the facts that $q>1$ and $1/q=1/p-\alpha/n$. Since $q/p>1$, by using the well-known inequality $\sum_i(\mu_i)^{q/p}\le(\sum_i\mu_i)^{q/p}$, we conclude that 
\begin{align}
J''_2&\le C\sum_{k=k_0+1}^\infty2^{kq}\Big(\tau^{{p(k-k_0)}/n}\Big)^{n(1-q)}
\cdot\Big(\sum_i\big|Q^k_i\big|\Big)^{q/p}\notag\\
&\le C\sum_{k=k_0+1}^\infty2^{kq}\Big(\tau^{{p(k-k_0)}/n}\Big)^{n(1-q)}
\Big(2^{-kp}\big\|f\big\|^{p}_{WH^p}\Big)^{q/p}\notag\\
&\le C\big\|f\big\|^{q}_{WH^p}\sum_{k=k_0+1}^\infty\Big(\tau^{{p(k-k_0)}/n}\Big)^{n(1-q)}\notag\\
&\le C\big\|f\big\|^{q}_{WH^p}.
\end{align}
Collecting the above inequality (4.5) with (4.3)--(4.4) and taking the supremum over all $\lambda>0$, and then taking $q$-th root on both sides, we finish the proof of Theorem 1.2.
\end{proof}

\section{Proof of Theorem 1.3}

\begin{proof}[Proof of Theorem 1.3]
Arguing as in the proof of Theorem 1.2, for any fixed $\lambda>0$, we can choose $k_0\in\mathbb Z$ satisfying $2^{k_0}\le\xi<2^{k_0+1}$, where we define $\xi=\lambda^{q/p}\big\|f\big\|^{1-q/p}_{WH^p}$. For every $f\in WH^p(\mathbb R^n)$, then in view of Theorem 2.2, we may write
\begin{equation*}
f=\sum_{k=-\infty}^\infty f_k=\sum_{k=-\infty}^{k_0} f_k+\sum_{k=k_0+1}^\infty f_k:=F_1+F_2,
\end{equation*}
where $F_1=\sum_{k=-\infty}^{k_0} f_k=\sum_{k=-\infty}^{k_0}\sum_i b^k_i$, $F_2=\sum_{k=k_0+1}^\infty f_k=\sum_{k=k_0+1}^\infty\sum_i b^k_i$ and $\{b^k_i\}$ satisfies $(a)$--$(c)$ in Theorem 2.2. Then we have
\begin{equation*}
\begin{split}
&\lambda^q\cdot \big|\big\{x\in\mathbb R^n:\big|T_{\Omega,\alpha}(f)(x)\big|>\lambda\big\}\big|\\
\le\,&\lambda^q\cdot \big|\big\{x\in\mathbb R^n:\big|T_{\Omega,\alpha}(F_1)(x)\big|>\lambda/2\big\}\big|
+\lambda^q\cdot \big|\big\{x\in\mathbb R^n:\big|T_{\Omega,\alpha}(F_2)(x)\big|>\lambda/2\big\}\big|\\
=\,&K_1+K_2.
\end{split}
\end{equation*}
Since $\Omega\in\mbox{Din}^r_\beta(S^{n-1})$ with $r>n/{(n-\alpha)}$, this is equivalent to $1\le r'<n/{\alpha}$. Then we are able to find a positive number $p_1$ such that $1\le r'<p_1<n/{\alpha}$. We also take $q_1>p_1>1$ such that $1/{q_1}=1/{p_1}-\alpha/n$. Hence, by Theorem B, we obtain that $T_{\Omega,\alpha}$ is bounded from $L^{p_1}(\mathbb R^n)$ to $L^{q_1}(\mathbb R^n)$. Repeating the arguments used in the proof of Theorem 1.2, we can also show that
\begin{equation*}
K_1\le C\big\|f\big\|^{q}_{WH^p}.
\end{equation*}
Let us now consider the other term $K_2$. As before, we set
\begin{equation*}
A_{k_0}=\bigcup_{k=k_0+1}^\infty\bigcup_i \widetilde{Q^k_i},
\end{equation*}
where $\widetilde{Q^k_i}=Q\big(x^k_i,\tau^{{p(k-k_0)}/n}(2\sqrt n)r^k_i\big)$ and $\tau$ is an appropriately chosen number such that $1<\tau<2$. Thus, we can further decompose $K_2$ as
\begin{equation*}
\begin{split}
K_2&\le\lambda^q\cdot \big|\big\{x\in A_{k_0}:|T_{\Omega,\alpha}(F_2)(x)|>\lambda/2\big\}\big|+
\lambda^q\cdot \big|\big\{x\in (A_{k_0})^c:|T_{\Omega,\alpha}(F_2)(x)|>\lambda/2\big\}\big|\\
&=K'_2+K''_2.
\end{split}
\end{equation*}
By using the same procedure as in Theorem 1.2, we can also obtain
\begin{equation*}
K'_2\le C\big\|f\big\|^{q}_{WH^p}.
\end{equation*}
It remains to estimate the last term $K''_2$. We first apply Chebyshev's inequality to obtain 
\begin{equation*}
\begin{split}
K''_2&\le 2^q\int_{(A_{k_0})^c}\big|T_{\Omega,\alpha}(F_2)(x)\big|^q\,dx\\
&\le 2^q\sum_{k=k_0+1}^\infty\sum_i
\int_{\big(\widetilde{Q^k_i}\big)^c}\big|T_{\Omega,\alpha}\big(b^k_i\big)(x)\big|^q\,dx\\
&= 2^q\sum_{k=k_0+1}^\infty\sum_i {\mathcal J}^k_i.
\end{split}
\end{equation*}
Again we denote $\tau^k_{i,\ell}=2^{\ell-1}\tau^{{p(k-k_0)}/n}\sqrt{n}r^k_i$ and
$$E^k_{i,\ell}=\big\{x\in\mathbb R^n:\tau^k_{i,\ell}\le|x-x^k_i|<2\tau^k_{i,\ell}\big\},\quad \ell=1,2,\ldots.$$
An application of H\"older's inequality leads to that
\begin{equation*}
\begin{split}
{\mathcal J}^k_i&\le\sum_{\ell=1}^\infty\int_{E^k_{i,\ell}}\big|T_{\Omega,\alpha}\big(b^k_i\big)(x)\big|^q\,dx\\
&\le\sum_{\ell=1}^\infty\Big|E^k_{i,\ell}\Big|^{1-q}
\bigg(\int_{E^k_{i,\ell}}\big|T_{\Omega,\alpha}\big(b^k_i\big)(x)\big|\,dx\bigg)^q.
\end{split}
\end{equation*}
Notice that $[n(1/p-1)]=0$ by the hypothesis. Consequently, by the cancellation condition of $b^k_i\in L^\infty(\mathbb R^n)$, we can get
\begin{equation*}
\begin{split}
\int_{E^k_{i,\ell}}\big|T_{\Omega,\alpha}\big(b^k_i\big)(x)\big|\,dx
&=\int_{E^k_{i,\ell}}\left|\int_{Q^k_i}
\Big[K_\alpha\big(x,x-y\big)-K_\alpha\big(x,x-x^k_i\big)\Big]b^k_i(y)\,dy\right|\,dx\\
&\le\int_{Q^k_i}\bigg\{\int_{E^k_{i,\ell}}\Big|K_\alpha\big(x,x-y\big)-K_\alpha\big(x,x-x^k_i\big)\Big|\,dx\bigg\}
\big|b^k_i(y)\big|\,dy\\
&\le\big\|b^k_i\big\|_{L^\infty}\int_{Q^k_i}
\bigg\{\int_{E^k_{i,\ell}}\Big|K_\alpha\big(x,x-y\big)-K_\alpha\big(x,x-x^k_i\big)\Big|\,dx\bigg\}\,dy.
\end{split}
\end{equation*}
If $y\in Q^k_i$ and $x\in\big(\widetilde{Q^k_i}\big)^c$, then we still have $\big|x-x^k_i\big|\ge 2\big|y-x^k_i\big|$ for all $i$ and $k$ by (3.4). Applying H\"older's inequality with exponent $r>1$ and Lemma 2.1, we can see that for any $y\in Q^k_i$, the integral in the brace of the above expression is bounded by
\begin{equation*}
\begin{split}
&\bigg(\int_{E^k_{i,\ell}}\Big|K_\alpha\big(x,x-y\big)-K_\alpha\big(x,x-x^k_i\big)\Big|^r\,dx\bigg)^{1/r}
\bigg(\int_{E^k_{i,\ell}}1\,dx\bigg)^{1/{r'}}\\
\le\, &C\cdot\Big|E^k_{i,\ell}\Big|^{1/{r'}}
\left(\int_{\tau^k_{i,\ell}\le|x|<2\tau^k_{i,\ell}}\Big|K_\alpha\big(x+x^k_i,x-(y-x^k_i)\big)
-K_\alpha\big(x+x^k_i,x\big)\Big|^r\,dx\right)^{1/r}\\
\le\, &C\cdot\Big|E^k_{i,\ell}\Big|^{1/{r'}}
\cdot\Big(\tau^k_{i,\ell}\Big)^{n/r-(n-\alpha)}
\left(\frac{|y-x^k_i|}{\tau^k_{i,\ell}}+
\int_{|y-x^k_i|/{2\tau^k_{i,\ell}}}^{|y-x^k_i|/{\tau^k_{i,\ell}}}\frac{\omega_r(\delta)}{\delta}\,d\delta\right)
\\
\le\, &C\cdot\Big|E^k_{i,\ell}\Big|^{1/{r'}}
\cdot\Big(\tau^k_{i,\ell}\Big)^{n/r-(n-\alpha)}
\left(\frac{|y-x^k_i|}{\tau^k_{i,\ell}}+\frac{|y-x^k_i|^\beta}{(\tau^k_{i,\ell})^\beta}\times
\int_{|y-x^k_i|/{2\tau^k_{i,\ell}}}^{|y-x^k_i|/{\tau^k_{i,\ell}}}
\frac{\omega_r(\delta)}{\delta^{1+\beta}}\,d\delta\right)\\
\le\, & C\cdot\Big(2\tau^k_{i,\ell}\Big)^{n/{r'}}
\cdot\Big(\tau^k_{i,\ell}\Big)^{n/r-(n-\alpha)}
\left(\frac{1}{2^{\ell}\tau^{{p(k-k_0)}/n}}+
\Big[\frac{1}{2^{\ell}\tau^{{p(k-k_0)}/n}}\Big]^\beta
\int_{0}^{1}\frac{\omega_r(\delta)}{\delta^{1+\beta}}\,d\delta\right)\\
\le\, & C\cdot\Big(2\tau^k_{i,\ell}\Big)^{\alpha}
\left(1+\int_{0}^{1}\frac{\omega_r(\delta)}{\delta^{1+\beta}}\,d\delta\right)
\times\left(\frac{1}{2^{\ell}\tau^{{p(k-k_0)}/n}}\right)^\beta.
\end{split}
\end{equation*}
Recall that $\big\|b^k_i\big\|_{L^\infty}\le C 2^k$ and $q\le1$. From the above estimate, it follows that
\begin{equation*}
\begin{split}
{\mathcal J}^k_i&\le C\cdot2^{kq}\sum_{\ell=1}^\infty \Big|E^k_{i,\ell}\Big|^{1-q}
\Big(2\tau^k_{i,\ell}\Big)^{\alpha q}\cdot\Big|Q^k_i\Big|^q
\left(\frac{1}{2^{\ell}\tau^{{p(k-k_0)}/n}}\right)^{\beta q}\\
&\le C\cdot2^{kq}\sum_{\ell=1}^\infty
\Big(2\tau^k_{i,\ell}\Big)^{n(1-q)+\alpha q}\cdot\Big|Q^k_i\Big|^q
\left(\frac{1}{2^{\ell}\tau^{{p(k-k_0)}/n}}\right)^{\beta q}\\
&\le C\cdot2^{kq}\sum_{\ell=1}^\infty \Big|Q^k_i\Big|^{1-q+{\alpha q}/n}\cdot\Big|Q^k_i\Big|^q
\Big(2^{\ell}\tau^{{p(k-k_0)}/n}\Big)^{n(1-q)+\alpha q-\beta q}\\
&\le C\cdot2^{kq}\cdot\Big|Q^k_i\Big|^{1+{\alpha q}/n}\Big(\tau^{{p(k-k_0)}/n}\Big)^{n(1-q)+\alpha q-\beta q},
\end{split}
\end{equation*}
where the last inequality holds since $q>n/{(n-\alpha+\beta)}$. Using the fact that $1+{\alpha q}/n=q/p>1$ and the well-known inequality $\sum_i(\mu_i)^{q/p}\le(\sum_i\mu_i)^{q/p}$, we deduce that
\begin{equation*}
\begin{split}
K''_2&\le C\sum_{k=k_0+1}^\infty\sum_i 2^{kq}\cdot\Big|Q^k_i\Big|^{1+{\alpha q}/n}
\Big(\tau^{{p(k-k_0)}/n}\Big)^{n-(n-\alpha+\beta)q}\\
&\le C\sum_{k=k_0+1}^\infty 2^{kq}\Big(\tau^{{p(k-k_0)}/n}\Big)^{n-(n-\alpha+\beta)q}
\cdot\Big(\sum_i\big|Q^k_i\big|\Big)^{q/p}\\
&\le C\big\|f\big\|^{q}_{WH^p}\sum_{k=k_0+1}^\infty\Big(\tau^{{p(k-k_0)}/n}\Big)^{n-(n-\alpha+\beta)q}\\
&\le C\big\|f\big\|^{q}_{WH^p}.
\end{split}
\end{equation*}
Summing up all the above estimates and taking the supremum over all $\lambda>0$, and then taking $q$-th root on both sides, we conclude the proof of Theorem 1.3.
\end{proof}

We finally remark that for any function $f$, a straightforward computation shows that the grand maximal function of $f$ is pointwise dominated by $M(f)$, where $M$ denotes the standard Hardy--Littlewood maximal operator. Hence, by the weak type (1,1) estimate of $M$, it is easy to see that the space $L^1(\mathbb R^n)$ is continuously embedded as a subspace of $WH^1(\mathbb R^n)$, and we have $\|f\|_{WH^1}\le C\|f\|_{L^1}$ for any $f\in L^1(\mathbb R^n)$. Therefore, as direct consequences of Theorems 1.1 and 1.2, we immediately obtain the following result.

\newtheorem{corollary}[theorem]{Corollary}

\begin{corollary}
Let $0<\beta\le 1$ and $p=1$. Suppose that $\Omega\in\mbox{Din}^r_\beta(S^{n-1})$ with $r>{2(n-1)}/n$, then there exists a constant $C>0$ independent of $f$ such that
\begin{equation*}
\big\|T_\Omega(f)\big\|_{WL^1}\le C\big\|f\big\|_{L^1}.
\end{equation*}
\end{corollary}

\begin{corollary}
Let $0<\alpha\le 1$, $p=1$ and $1/q=1/p-\alpha/n$. Suppose that $\Omega\in\mbox{Din}^{n/{(n-\alpha)}}_\alpha(S^{n-1})$, then there exists a constant $C>0$ independent of $f$ such that
\begin{equation*}
\big\|T_{\Omega,\alpha}(f)\big\|_{WL^{n/{(n-\alpha)}}}\le C\big\|f\big\|_{L^1}.
\end{equation*}
\end{corollary}


\begin{thebibliography}{99}

\bibitem{cal1} A. P. Calder\'on and A. Zygmund, On a problem of Mihlin, Trans. Amer. Math. Soc, \textbf{78}(1955), 209--224.
\bibitem{cal2} A. P. Calder\'on and A. Zygmund, On singular integrals, Amer. J. Math, \textbf{78}(1956), 289--309.
\bibitem{cal3} A. P. Calder\'on and A. Zygmund, On singular integrals with variable kernels, Appl. Anal, \textbf{7}(1978), 221--238.
\bibitem{ding6} Y. Ding, C. C. Lin and Y. C. Lin, Erratum: ``On Marcinkiewicz integral with variable kernels, Indiana Univ. Math. J, \textbf{53}(2004), 805--821", Indiana Univ. Math. J, \textbf{56}(2007), 991--994.
\bibitem{ding7} Y. Ding, C. C. Lin and S. L. Shao, On the Marcinkiewicz integral with variable kernels, Indiana Univ. Math. J, \textbf{53}(2004), 805--821.
\bibitem{ding1} Y. Ding, S. Z. Lu and S. L. Shao, Integral operators with variable kernels on weak Hardy spaces, J. Math. Anal. Appl, \textbf{317}(2006), 127--135.
\bibitem{ding2} Y. Ding, S. Z. Lu and Q. Y. Xue, Marcinkiewicz integral on Hardy spaces, Integr. Equ. Oper. Theory, \textbf{42}(2002), 174--182.
\bibitem{ding3} Y. Ding, S. Z. Lu and Q. Y. Xue, Parametrized Littlewood--Paley operators on Hardy and weak Hardy spaces, Math. Nachr, \textbf{280}(2007), 351--363.
\bibitem{ding4} Y. Ding, S. Z. Lu and Q. Y. Xue, Parametrized area integrals on Hardy spaces and weak Hardy spaces, Acta Math. Sinica (Engl. Ser.), \textbf{23}(2007), 1537--1552.
\bibitem{ding5} Y. Ding and X. F. Wu, Weak Hardy space and endpoint estimates for singular integrals on space of homogeneous type, Turkish J. Math, \textbf{34}(2010), 235--247.
\bibitem{cfefferman} C. Fefferman, N. Rivi\`ere and Y. Sagher, Interpolation between $H^p$ spaces: The real method, Trans. Amer. Math. Soc, \textbf{191}(1974), 75--81.
\bibitem{rfefferman} R. Fefferman and F. Soria, The space weak $H^1$, Studia Math, \textbf{85}(1987), 1--16.
\bibitem{kurtz} D. S. Kurtz and R. L. Wheeden, Results on weighted norm inequalities for multipliers, Trans. Amer. Math. Soc, \textbf{255}(1979), 343--362.
\bibitem{lee} M. Y. Lee, C. C. Lin, Y. C. Lin and D. Y. Yan, Boundedness of singular integral operators with variable kernels, J. Math. Anal. Appl, \textbf{348}(2008), 787--796.
\bibitem{lin} C. C. Lin, Y. C. Lin, X. X. Tao and X. Yu, The boundedness of Marcinkiewicz integral with variable kernel, Illinois J. Math, \textbf{53}(2009), 197--217.
\bibitem{liu1} H. P. Liu, The weak $H^p$ spaces on homogeneous groups, Lecture Notes in Math, Vol. 1494, Springer-Verlag, 1991, 113--118.
\bibitem{liu2} H. P. Liu, The wavelet characterization of the space Weak $H^1$, Studia Math, \textbf{103}(1992), 109--117.
\bibitem{lu} S. Z. Lu, Four Lectures on Real $H^p$ Spaces, World Scientific Publishing, River Edge, N.J., 1995.
\bibitem{muckenhoupt} B. Muckenhoupt and R. L. Wheeden, Weighted norm inequalities for singular and fractional integrals, Trans. Amer. Math. Soc, \textbf{161}(1971), 249--258.
\bibitem{tao} X. X. Tao, X. Yu and S. Y. Zhang, Marcinkiewicz integrals with variable kernels on Hardy and weak Hardy spaces, J. Funct. Spaces Appl, \textbf{8}(2010), 1--16.

\end{thebibliography}
\end{document}